\documentclass[11pt]{amsart}
\theoremstyle{plain}

\newtheorem{thm}{Theorem}[section]
\newtheorem{theorem}[thm]{Theorem}

\newtheorem{corollary}[thm]{Corollary}
\newtheorem{lemma}[thm]{Lemma}

\newtheorem{prop}[thm]{Proposition}
\newtheorem{proposition}[thm]{Proposition}

\newtheorem*{conjecture*}{Conjecture}

\newtheorem*{question*}{Question}

\theoremstyle{definition}
\newtheorem{defn}[thm]{Definition}

\newtheorem*{rem*}{Remark}
\newtheorem{remark}[thm]{Remark}
\newtheorem*{remark*}{Remark}
\newtheorem{remarks}[thm]{Remarks}
\newtheorem*{remarks*}{Remarks}

\newtheorem*{example*}{Example}

\newtheorem*{examples*}{Examples}

\newtheorem*{notation*}{Notation}

\newtheorem*{bibliographical-note}{Bibliographical note}
\newcommand{\zeroindent}{\parindent0cm \parskip1ex}
\newcommand{\acknowledgements}{{\em Acknowledgements. }}

{\begin{list}{(\alph{enumi}) }{\usecounter{enumi}

\leftmargin0cm \labelsep0cm \rightmargin0cm
\setlength{\labelwidth}{\fill}}}
{\end{list}}

\newenvironment{normallist}%
{

\begin{enumerate}}{\end{enumerate}}

{

\begin{enumerate}}{\end{enumerate}}

\newenvironment{romanlist}%
{

\begin{enumerate} \parsep0cm
\itemsep0cm \parskip0cm}{\end{enumerate}}

\newenvironment{primelist}%
{

\begin{enumerate}}{\end{enumerate}}

{

\begin{enumerate}}{\end{enumerate}}

\newcommand{\R}{\mathbb{R}}
\newcommand{\Z}{\mathbb{Z}}

\newcommand{\C}{\mathbb{C}}
\newcommand{\N}{\mathbb{N}}

\newcommand{\id}{\mathrm{id}}
\newcommand{\im}{\mathrm{im}}

\newcommand{\suchthat}{\; | \;}
\newcommand{\iso}{\cong}           
\newcommand{\smooth}{C^\infty}
\newcommand{\CP}[1]{\C {\mathrm P}^{#1}}

\newcommand{\half}{{\textstyle\frac{1}{2}}}
\newcommand{\quarter}{{\textstyle\frac{1}{4}}}
\newcommand{\mymod}{\quad\text{mod }}
\newcommand{\Rgeq}{\R^{\scriptscriptstyle \geq 0}}

\newcommand{\mo}{(M,\omega)}
\renewcommand{\o}{\omega}

\newcommand{\Aut}{\mathrm{Aut}}

\newcommand{\Symp}{\mathrm{Sp}}
\newcommand{\leftsc}{\langle}
\newcommand{\rightsc}{\rangle}
\newcommand{\Diff}{\mathrm{Diff}}

\zeroindent
\numberwithin{equation}{section}

\theoremstyle{plain}

\newenvironment{remarkslist}%
{\begin{list}{(\alph{enumi}) }{\usecounter{enumi}

\leftmargin0cm \labelsep0cm
\itemsep1em \rightmargin0cm
\setlength{\labelwidth}{\fill}}}
{\end{list}}

\newcommand{\tL}{\widetilde{L}}
\renewcommand{\P}{\mathcal{P}}
\newcommand{\PP}{\P(L,L')}

\renewcommand{\H}{\mathcal H}
\newcommand{\reg}{\mathrm{reg}}
\newcommand{\JJ}{\mathcal{J}}
\newcommand{\J}{\mathbf{J}}
\newcommand{\moduli}{\mathcal{M}}
\newcommand{\V}{\mathcal V}
\newcommand{\loc}{\mathrm{loc}}
\newcommand{\cotangent}{T^{*\!}}
\newcommand{\disc}{D_{\epsilon}(\cotangent S^2)}
\renewcommand{\L}{\mathcal{L}}

\newcommand{\surgery}{\#}

\newcommand{\One}{\mathbf{1}}
\newcommand{\gen}[1]{\leftsc #1 \rightsc}

\title[Lagrangian two-spheres]{Lagrangian two-spheres
can be symplectically knotted}
\author{Paul Seidel}
\date{19/3/98}
\thanks{Supported by the European Community and by NSF grant DMS 9304580.}
\address{Centre de Math{\'e}matiques, Ecole Polytechnique, Palaiseau, France}
\email{seidel@math.polytechnique.fr}
\begin{document}
\maketitle
\section{Introduction}

In the past few years there have been several striking results about the topology of
Lagrangian surfaces in symplectic four-manifolds. The general tendency of these results is
that many isotopy classes of embedded surfaces do not contain Lagrangian representatives. This
is called the {\em topological unknottedness} of Lagrangian surfaces; see
\cite{eliashberg-polterovich93b} for a survey. The aim of this paper is to complement this
picture by showing that Lagrangian surfaces can be {\em symplectically knotted} in infinitely
many inequivalent ways. That is to say, a single isotopy class of embedded surfaces can
contain infinitely many Lagrangian representatives which are non-isotopic in the Lagrangian
sense. The symplectic four-manifolds for which we prove this are non-compact in a mild sense;
they are interiors of compact symplectic manifolds with contact type boundary. For instance,
one can take
\begin{equation} \label{eq:affine}
M = \{ z \in \C^3 \suchthat |z| < R, \quad
z_1^2 + z_2^2 = z_3^{m+1} + \half\}
\end{equation}
for $m \geq 3$ and large $R$, with the standard symplectic form. It
seems likely that the same phenomenon occurs for a large class of
closed symplectic four-manifolds. At present, the best result in this
direction is that for any $N$, there is a closed four-manifold
containing $N$ Lagrangian two-spheres which
are all isotopic as smooth submanifolds but pairwise non-isotopic as
Lagrangian submanifolds. For example, the smooth hypersurface in
$\CP{3}$ of degree $N+4$ has this property. We will not prove the
statement about closed manifolds in this paper.

As a by-product of our construction it turns out that the
four-manifolds which we consider have a symplectic automorphism
$\phi$ with the the following property: $\phi$ is isotopic to the identity
as a diffeomorphism, but not symplectically so. Moreover, none of
the iterates $\phi^r$ are symplectically isotopic to the identity.
To be precise, $\phi$ is the identity outside a compact subset,
and it can be deformed to the identity map inside the group of
diffeomorphisms with this property. In contrast, it cannot be
deformed to the identity through symplectomorphisms even if we allow
these to behave arbitrarily at infinity; and the same holds for
$\phi^r$. Just as in the case of Lagrangian two-spheres, one can obtain a
weaker form of this statement for closed symplectic four-manifolds;
see \cite{seidel97}. Kronheimer has obtained results of a related
kind using a parametrized version of Seiberg-Witten theory \cite{kronheimer97}.

We recall some basic definitions. Let $M$ be a
differentiable manifold. A differentiable isotopy between two compact
submanifolds $L_0,L_1 \subset M$ is a compact submanifold
$\tL \subset M \times [0;1]$ with $\tL \cap (M \times \{t\}) = L_t
\times \{t\}$  for $t = 0,1$ and such that the intersection
$\tL \cap (M \times \{t\})$ is transverse for all $t \in [0;1]$.
If $\mo$ is a symplectic manifold, a Lagrangian isotopy between
two compact Lagrangian submanifolds $L_0,L_1 \subset M$ is an
$\tL$ as above, with the additional property that $\tL \cap
(M \times \{t\}) \subset M \times \{t\}$ is Lagrangian for all $t$.
This means that the pullback of $\o$ to $\tL$ via the projection
$\tL \subset M \times [0;1] \longrightarrow M$ is of the form
$\theta \wedge dt$ for some $\theta  \in \Omega^1(\tL)$. The
Lagrangian isotopy $\tL$ is called exact if one can write
$\theta \wedge dt = d(H dt)$ for some $H \in \smooth(\tL,\R)$.
Any Lagrangian isotopy between Lagrangian submanifolds with
vanishing first Betti number is exact.

Two diffeomorphisms $\phi_0,\phi_1: M \longrightarrow M$ are differentiably isotopic if they
can be joined by a path $(\phi_t)_{0 \leq t \leq 1}$ in $\Diff(M)$ which is smooth (in the
sense that $F(x,t) = \phi_t(x)$ is a smooth map from $M \times [0;1]$ to $M$). Two compact
submanifolds $L_0,L_1$ which are differentiably isotopic are also ambient isotopic, that is,
there is a path $(\phi_t)$ with $\phi_0 = \id$ and $\phi_1(L_0) = L_1$. Similarly, one can
define the notion of symplectic isotopy between two symplectic automorphisms, and two compact
Lagrangian submanifolds which are exact Lagrangian isotopic are also ambient isotopic in the
symplectic sense.

We now explain how our examples of symplectically knotted Lagrangian
two-spheres are constructed. Given a Lagrangian
two-sphere $L$ in a symplectic four-manifold $\mo$, one can
define a symplectic automorphism $\tau_L$ of $M$ called the
{\em generalized Dehn twist along $L$}. This automorphism is
trivial outside a tubular neighbourhood of $L$, and its
restriction to $L$ is the antipodal involution of $S^2$.
The definition of $\tau_L$ involves additional choices, but
the outcome is independent of these choices up to symplectic
isotopy. The topological analogues of generalized Dehn twists
are certain diffeomorphisms associated to embedded two-spheres of
self-intersection $-2$ in smooth four-manifolds. These maps
are familiar to topologists. The symplectic viewpoint appears
for the first time in Arnol'd's paper \cite{arnold95}.

We will use these generalized Dehn twists in the following way:
assume that $M$ contains two Lagrangian two-spheres $L_1,L_2$.
Then one can construct an infinite family of such two-spheres
\begin{equation} \label{eq:loner}
L_1^{(r)} = \tau_{L_2}^{2r}(L_1) \quad (r \in \Z)
\end{equation}
by twisting $L_1$ around $L_2$. We have used only even iterates
of $\tau_{L_2}$ because the square of any generalized Dehn twist is
differentiably isotopic to the identity (this was explained to the
author by Peter Kronheimer). As a consequence all the
$L_1^{(r)}$ are isotopic as differentiable submanifolds. Our
main result is that they are not always Lagrangian isotopic.

\begin{theorem}\label{th:main}
Let $\mo$ be the interior of a compact symplectic manifold with
contact type boundary. We assume that $\mo$ contains three Lagrangian
two-spheres $L_1,L_2,L_3$ such that $L_1 \cap L_3 = \emptyset$,
and such that both $L_1 \cap L_2$ and $L_2 \cap L_3$ are transverse
and consist of a single point (this is called an $(A_3)$-configuration
of Lagrangian two-spheres). Moreover, we require that both
$[\o] \in H^2(M;\R)$ and $c_1\mo \in H^2(M;\Z)$ vanish.
Let $L_1^{(r)}$ be as in \eqref{eq:loner}. Then $L_1^{(s)}$ and
$L_1^{(t)}$ are not Lagrangian isotopic unless $s = t$. \end{theorem}

The assumptions are very restrictive. To provide some concrete
examples, we will prove that the affine hypersurfaces
\eqref{eq:affine} satisfy them. As an immediate consequence,
one obtains the following result which establishes our
second claim:

\begin{corollary} \label{th:corollary}
In the situation of Theorem \ref{th:main} no iterate of $\tau_{L_2}^2$
is symplectically isotopic to the identity. \qed \end{corollary}

\acknowledgements I am indebted to Peter Kronheimer for several
helpful conversations, and to Leonid Polterovich for introducing me
to Lagrangian surgery. Parts of this paper were written during
visits to Tel Aviv University and Stanford University.

\section{An outline of the proof \label{sec:outline}}

Let $\mo$ be the interior of a compact symplectic four-manifold with
contact type boundary, and assume that $[\o] \in H^2(M;\R)$ is zero.
The Floer homology $HF(L,L')$ is a finite-dimensional $\Z/2$-vector
space associated to any pair $(L,L')$ of Lagrangian two-spheres in
$M$. It is invariant under Lagrangian isotopy in the following sense:
if $L''$ is Lagrangian isotopic to $L$ then $HF(L,L') \iso
HF(L'',L')$ for all $L'$. Conversely, to prove that two Lagrangian
two-spheres $L$ and $L''$ are not Lagrangian isotopic, it is
sufficient to find a third two-sphere $L'$ such that $HF(L,L')
\not\iso HF(L'',L')$. The main difficulty is that the Floer homology
groups are difficult to compute. If $L$ and $L'$ are disjoint,
$HF(L,L') = 0$. Floer \cite{floer89} proved that if $L$ and $L'$
are Lagrangian isotopic then $HF(L,L') \iso H_*(L;\Z/2)$. We will
use a generalization of Floer's result due to Pozniak \cite{pozniak}.

\begin{defn} $L,L' \subset M$ have {\em clean intersection} if $
N = L \cap L'$ is a smooth submanifold of $M$ and satisfies
$TN = (TL|N) \cap (TL'|N)$.
\end{defn}

To simplify the statement, we consider only the case when (in
addition to the conditions imposed above) the first Chern class
of $M$ vanishes. Then the Floer homology of $(L,L')$ is a graded
group. The grading is not quite unique. However, this is
irrelevant for our purpose since we use it only as a computational
device: ultimately, only the ungraded Floer homology group will serve as an
invariant. Pozniak's result can be formulated as follows: if $L$
and $L'$ have clean intersection, there is a spectral sequence which
converges to $HF_*(L,L')$ and whose $E^1$-term is
\[
E^1_{pq} = \begin{cases}
H_{p+q-i'(C_p)}(C_p;\Z/2) & 1 \leq p \leq r,\\
0 & \text{otherwise.}
\end{cases}
\]
Here $C_1, \dots, C_r$ are the connected components of $L \cap L'$
ordered in a way determined by the action functional. $i'(C_p) \in \Z$
is a kind of Maslov index. This is a homology spectral sequence, that is,
the $d$-th differential ($d \geq 1$) has degree $(-d,d-1)$.

In the situation of Theorem \ref{th:main} one can arrange that for
a given $r>0$, $L_1^{(r)}$ intersects $L_3$ cleanly in $r$ circles
$C_1,\dots, C_r$. With the right ordering one finds that $i'(C_p) = 2p$.
Hence the $E^1$ term of the Pozniak spectral sequence for
$(L_1^{(r)},L_3)$ is
\[
E^1_{pq} = \begin{cases}
\Z/2 & 1 \leq p \leq r \text{ and } q = p \text{ or } p+1,\\
0 & \text{otherwise.}
\end{cases}
\]
Since the entry $E^1_{11} = \Z/2$ survives to $E^{\infty}$, we have
$HF(L_1^{(r)},L_3) \neq 0$. On the other hand, $HF(L_1,L_3) = 0$
because $L_1$ and $L_3$ are assumed to be disjoint. Hence $L_1$ is not Lagrangian
isotopic to $L_1^{(r)}$, for any $r>0$. It follows that $L_1^{(s)}$ is
not Lagrangian isotopic to $L_1^{(t)}$ whenever $s<t$ because otherwise
$L_1^{(t-s)}$ would be Lagrangian isotopic to $L_1$.

\section{Floer homology and clean intersection\label{sec:pozniak}}

This section and the next two contain a more detailed account
of Floer homology and of Pozniak's results. Our definition of
Floer homology is essentially Floer's original one \cite{floer88c}.
The construction has been generalized by Oh \cite{oh93} but this
generalization is unnecessary for our purpose.

Let $(\overline{M},\bar{\o})$ be a compact symplectic manifold with
contact type boundary, and $\mo$ its interior. Fix an $\bar{\o}$-compatible
almost complex structure $\overline{J}$ which makes the boundary
$\overline{J}$-convex.

Let $L,L' \subset M$
be a pair of compact Lagrangian submanifolds, and $\PP$ the space
of smooth paths $\gamma: I = [0;1] \longrightarrow M$ such that
$\gamma(0) \in L$ and $\gamma(1) \in L'$. Let
\[
\alpha(\gamma)\xi = \int_I \o(\frac{d\gamma}{dt},\xi(t))\, dt
\]
be the action one-form on $\PP$. $\alpha$ is always closed.
The Floer homology group $HF(L,L')$ is defined whenever $[\alpha]
\in H^1(\PP;\R)$ is zero. If this holds for $(L,L')$, it also holds
for $(L'',L')$ whenever $L''$ is exact Lagrangian isotopic to $L$,
and then $HF(L'',L') \iso HF(L,L')$. More precisely, an
exact Lagrangian isotopy from $L$ to $L''$ determines an isomorphism
between the two Floer homologies. If we assume that $[\o] = 0$ and
$H^1(L;\R) = H^1(L';\R) = 0$, then $[\alpha] = 0$ and all Lagrangian
isotopies are exact. In this way one recovers the description of
Floer homology given in the previous section.

We will now review briefly the definition of $HF(L,L')$. By
assumption, one can choose a function $a: \PP \longrightarrow \R$
with $da = \alpha$. For any $H \in \H = \smooth_c(I \times M,\R)$, let
\[
a_H(\gamma) = a(\gamma) + \int_I H_t(\gamma(t))\, dt.
\]
Let $Z(H) \subset \PP$ be the set of critical points of $a_H$.
For $H = 0$ the critical points are the constant paths $\gamma_x$
at $x \in L \cap L'$. In general, $\gamma$ is a critical point iff
it is an orbit of the flow $(\phi_t^H)$ induced by $H$. Hence $Z(H)$
can be identified naturally with $\phi_1^H(L) \cap L'$.
Let $\H^{\reg} \subset \H$
be the dense subset of those $H$ for which $\phi_1^H(L)$ and
$L'$ intersect transversely. For $H \in \H^{\reg}$ one defines
$CF(H)$ to be the $\Z/2$-vector space generated by the finite
set $Z(H)$. Let $\JJ$ be the space of one-parameter families
$\J = (J_t)_{t \in I}$ of almost complex structures on $M$ which
have the following properties:

\begin{normallist}
\item Every $J_t$ is $\o$-compatible.
\item There is a neighbourhood $U \subset \overline{M}$ of
$\partial\overline{M}$ such that
$J_t|U \cap M = \overline{J}|U \cap M$ for all $t$.
\end{normallist}

For $H \in \H$, $\J \in \JJ$, and $\gamma_-,\gamma_+ \in Z(H)$, we
denote by $\moduli_{H,\J}(\gamma_-,\gamma_+)$ the set of smooth maps
$u: \R \times I \longrightarrow M$ such that
\begin{equation}
\label{eq:floer}
\begin{cases}
u(s,0) \in L, \quad u(s,1) \in L', \quad
\lim_{s \rightarrow \pm\infty} u(s,\cdot) = \gamma_\pm, &\\
\frac{\partial u}{\partial s} +
J_t(u) \left(\frac{\partial u}{\partial t} - X_t^H(u)\right) = 0.&
\end{cases}
\end{equation}
Here $X^H$ is the (time-dependent) Hamiltonian vector field of $H$. If one
thinks of $u$ as a map $\R \longrightarrow \PP$, the solutions
of \eqref{eq:floer} are the bounded negative gradient flow lines of
$a_H$ with respect to an $L^2$-metric defined by $\J$. Hence
$\moduli_{H,\J}(\gamma_-,\gamma_+)$ can be nonempty only if
$a_H(\gamma_-) > a_H(\gamma_+)$ or $\gamma_- = \gamma_+$.
The assumption on the behaviour of $\J$ at infinity implies that
the union of the images of all solutions of \eqref{eq:floer} lies
inside a compact subset $K \subset M$. This removes any possible problems
arising from the non-compactness of $M$.

Assume that $H \in \H^{\reg}$, and let $\JJ^{\reg}(H) \subset \JJ$
be the subspace of those $\J$ for which all solutions of \eqref{eq:floer}
are reguar. Regularity is defined as the surjectivity of the
linearization of \eqref{eq:floer} in suitable Sobolev spaces.
$\JJ^{\reg}(H)$ is a dense subset; the necessary transversality
arguments were carried out in \cite{floer-hofer-salamon94} and \cite{oh96c}.
If $\J$ is in $\JJ^{\reg}(H)$ the spaces $\moduli_{H,\J}(\gamma_-,\gamma_+)$
have a natural structure of finite-dimensional smooth manifolds.
Moreover, any one of them has only finitely many one-dimensional
connected components. Let $n_{H,\J}(\gamma_-,\gamma_+) \in \Z/2$ be
the number mod $2$ of these components. Floer proved that the
homomorphism
\[
\partial(H,\J): CF(H) \longrightarrow CF(H),
\quad \partial(H,\J) \gen{\gamma_-} = \sum_{\gamma_+ \in Z(H)}
n_{H,\J}(\gamma_-,\gamma_+) \gen{\gamma_+}
\]
satisfies $\partial(H,\J) \circ \partial(H,\J) = 0$. Floer
homology is defined by
\[
HF(L,L',H,\J) = \ker \partial(H,\J)/ \im \; \partial(H,\J).
\]
A continuation argument proves that this
is independent of the choice of $\J$ and $H$ up to canonical
isomorphisms. Therefore one can omit $(H,\J)$ from the notation.
By a simple change of variables, the independence of
$H$ implies the invariance under exact Lagrangian isotopy.
A detailed exposition of the continuation argument (in a slightly
different context) can be found in \cite{salamon-zehnder92}.

Floer \cite{floer88c} introduced versions of Floer homology which
are local with respect to certain parts of $N = L \cap L'$. More
precisely, let $C \subset N$ be a path component which is
both open and closed in $N$. Choose a $\J_0 \in \JJ$ and
an open neighbourhood
$U \subset M$ of $C$ such that $N \cap U = C$. There is
a contractible neighbourhood $\V \subset \H \times \JJ$ of $(0,\J_0)$
such that any $(H,\J) \in \V$ has the following properties:
\begin{normallist}
\item any $\gamma \in Z(H)$ such that $\im(\gamma) \subset
\overline{U}$ satisfies $\im(\gamma) \subset U$;
\item if $\gamma_-,\gamma_+ \in Z(H)$ satisfy $\im(\gamma_{\pm})
\subset U$ then any $u \in \moduli_{H,\J}(\gamma_-,\gamma_+)$
has $\im(u) \subset U$.
\end{normallist}
This can be proved by a simple limit argument. Now choose
$(H,\J) \in \V$ such that $H \in \H^{\reg}$ and $\J \in
\JJ^{\reg}(H)$. Let $CF^{\loc}(H;C)$ be the $\Z/2$-vector space
generated by those $\gamma \in Z(H)$ with $\im(\gamma) \subset U$,
and $\partial^{\loc}(H,\J;C): CF^{\loc}(H;C) \longrightarrow
CF^{\loc}(H;C)$ the homomorphism obtained by considering only
those $u \in \moduli_{H,\J}(\gamma_-,\gamma_+)$ such that
$\im(u) \subset U$. $\partial^{\loc}(H,\J;C)^2 = 0$, and one defines
$HF^{\loc}(L,L',H,\J;C) = \ker \partial^{\loc}(H,\J;C)/
\im \;\partial^{\loc}(H,\J;C)$.

One can prove that the local Floer homology is independent of
the choice of $(H,\J) \in \V$. The proof is again by a continuation
method. Recall that in arguments of this kind one studies maps
$u: \R \times I \longrightarrow M$ which satisfy an equation
\begin{equation} \label{eq:continuation}
\begin{cases}
u(s,0) \in L, \quad u(s,1) \in L', &\\
\frac{\partial u}{\partial s} +
J_{s,t}(u) \left(\frac{\partial u}{\partial t} - X_{s,t}(u)\right) = 0,&
\end{cases}
\end{equation}
with suitable asymptotic behaviour. Here $J_{s,t}$ is a two-parameter
family of almost complex structures and $X_{s,t}$ is the family of Hamiltonian
vector fields on $M$ induced by a function $K \in \smooth(\R \times I \times M,\R)$.
In the case of local Floer homology, one considers such families with the
property that $H_s = K_{s,\cdot} \in \H$ and $\J_s = J_{s,\cdot}$ satisfy
$(H_s,\J_s) \in \V$ for all $s$. In itself this does not imply that
any solution of \eqref{eq:continuation} with limits in $U$ satisfies $\im(u) \subset U$.
However, a limit argument shows that this holds if the path $s \longmapsto
(H_s,\J_s) \in \V$ is close to a constant path in a suitable sense. This is
sufficient to prove that the local Floer homology is independent of
$H$ and $\J$. The same argument also proves that it remains unchanged under
small variations of $\J_0$, and hence is independent of $\J_0$. Therefore one
obtains a well-defined group $HF^{\loc}(L,L';C)$. This
group is called local because it depends only on the behaviour
of $L$ and $L'$ in an arbitrarily small neighbourhood of $C$. For
instance, if $C = \{x\}$ is a transverse intersection point of $L$
and $L'$ then $HF^{\loc}(L,L';C) \iso \Z/2$.

\begin{theorem}[Pozniak \protect{\cite[Theorem 3.4.11]{pozniak}}]
\label{th:local-floer-homology} Assume that $L$ and $L'$ intersect
cleanly along $C$. Then $HF^{\loc}(L,L';C)$ is isomorphic to the
total homology $H_*(C;\Z/2)$. \end{theorem}

This is the main result of \cite{pozniak}. We will not reproduce
Pozniak's proof here. Instead, we will give (at the end of this
section) an indirect proof of a special case of Theorem
\ref{th:local-floer-homology}:

\begin{prop} \label{th:small-pozniak} Assume that $\mo$ is
four-dimensional, that $L$ and $L'$ are orientable, and that
$C$ is a circle along which $L$ and $L'$ intersect cleanly.
Then $HF^{\loc}(L,L';C) \iso \Z/2 \oplus \Z/2$. \end{prop}

The ordinary Floer homology and its local version are related in
the following way: let $L,L' \subset M$ be two compact Lagrangian
submanifolds such that $N = L \cap L'$ can be decomposed into finitely
many path components $C_1, \dots, C_r$, each of which is open
(and closed) in $N$. As always, we assume that $[\alpha]
\in H^1(\PP;\R)$ vanishes. Let $a_j = a(\gamma_{x_j})$ where
$\gamma_{x_j}$ is the constant path at a point $x_j \in C_j$.
We assume that the $C_j$ have been ordered in such a way that
$a_1 \leq a_2 \leq \dots \leq a_r$.
Choose a $\J_0 \in \JJ$ and open neighbourhoods $U_j \subset
M$ of $C_j$ whose closures are pairwise disjoint. There is a
contractible neighbourhood $\V \subset \H \times \JJ$ of $(0,\J_0)$
such that any $(H,\J) \in \V$ has the following properties:

\begin{normallist}
\item any $\gamma \in Z(H)$ satisfies $\im(\gamma) \subset U_j$ for
some $j$;
\item if $\gamma_-,\gamma_+ \in Z(H)$ satisfy $\im(\gamma_-) \subset
U_j$ and $\im(\gamma_+) \subset U_k$ with $j<k$ then
$\moduli_{H,\J}(\gamma_-,\gamma_+) = \emptyset$;
\item for $\gamma_-,\gamma_+ \in Z(H)$ such that $\im(\gamma_{\pm})
\subset U_j$ every $u \in \moduli_{H,\J}(\gamma_-,\gamma_+)$ satisfies
$\im(u) \subset U_j$.
\end{normallist}

Now take $(H,\J) \in \V$ such that $H \in \H^{\reg}$ and $\J \in
\JJ^{\reg}(H)$. Consider the filtration of $CF(H)$ by the subspaces
$CF(H)^{[j]}$ generated by those $\gamma \in Z(H)$ such that
$\im(\gamma) \subset U_1 \cup \dots \cup U_j$. It is a consequence of
the properties which we have just stated that $\partial(H,\J)$
preserves this filtration, and that the homology of the induced
boundary operator on $CF(H)^{[j]}/CF(H)^{[j-1]}$ is isomorphic to
the local Floer homology $HF^{\loc}(L,L';C_j)$.

\proof[Proof of Proposition \ref{th:small-pozniak}] The first
step in the proof is a local normal form theorem
\cite[Proposition 3.4.1]{pozniak} for cleanly intersecting
Lagrangian submanifolds. In our case this says that one can identify
a neighbourhood of $C$ in $M$ symplectically with a neighbourhood
of $S^1 \times 0$ in $(S^1 \times \R^3, ds \wedge dx_1
+ dx_2 \wedge dx_3)$ in such a way that $L$ is mapped to $S^1 \times
0 \times \R \times 0$ and $L'$ is mapped to $S^1 \times 0 \times 0
\times \R$. In particular, the local behaviour of $L$ and $L'$ near
$C$ is the same in all cases covered by Proposition
\ref{th:small-pozniak}. Hence the local
Floer homology group is also the same in all cases. We denote this
group, which we want to compute, by $G$.

The second step is to show that $G$ is either $0$ or $\Z/2 \oplus
\Z/2$. To do this, one observes that, given a neighbourhood $U \subset M$
of $C$, there are arbitrarily small $H \in \H$ such that there are
precisely two $\gamma \in Z(H)$ with $\im(\gamma) \in U$. In the local
model $S^1 \times \R^3$, one can take
\begin{equation} \label{eq:local-h}
H_t(z,x_1,x_2,x_3) = h(z) \psi(|x_1|^2 + |x_2|^2 + |x_3|^2),
\end{equation}
where $h$ is a Morse function on $S^1$ with two critical points
and $\psi$ is a cutoff function ($\psi(r) = 0$ for large $r$
and $= 1$ for small $r$).

The final step is to exclude the possibility that $G = 0$. Let
$M = \cotangent T^2$ with the standard symplectic structure. Take
a Morse-Bott function $k \in \smooth(T^2,\R)$ whose critical set
consists of two circles. Let $L \subset M$ be the zero-section
and $L' \subset M$ the graph of $dk$. Then $L'$ and $L$ intersect
cleanly in two circles. The considerations above show that for
suitable $H$ and $\J$, the chain group $CF(H)$ and the boundary
operator $\partial(H,\J)$ have the following form: there is a subgroup
$CF(H)^{[1]}$ which is preserved by $\partial(H,\J)$. The
homology of this subgroup is equal to $G$, and the homology of the
quotient is also equal to $G$. If we assume that $G$ is zero, the
long exact sequence would imply that $HF(L,L') = 0$. However,
since $L'$ is exact Lagrangian isotopic to $L$, Floer's theorem
says that $HF(L,L') \iso H_*(T^2;\Z_2)$. \qed

The proof of Proposition \ref{th:small-pozniak} was based on the
relationship between clean intersection in cotangent bundles and
Morse-Bott functions. In fact, one can see Pozniak's approach
as an analogue of the Morse-Bott spectral sequence in ordinary
homology (see \cite{austin-braam94} for an exposition). Other
arguments in Floer homology based on the same principles
can be found in \cite{ruan-tian95} and \cite{oh96}.

\section{The grading on Floer homology \label{sec:index}}

The material collected in this section is due to Viterbo
\cite{viterbo87}, Floer \cite{floer88e} and Robbin-Salamon
\cite{robbin-salamon93} \cite{robbin-salamon95}.

Let $\L(n)$ be the Lagrangian Grassmannian, which parametrizes
linear Lagrangian subspaces in $\R^{2n}$. The Maslov index
associates an integer $\mu(\lambda,\lambda')$ to a pair of loops
$\lambda,\lambda': S^1 \longrightarrow \L(n)$. This index
is invariant under homotopy and under conjugation of both
$\lambda$ and $\lambda'$ by a loop in $\Symp(2n;\R)$. Usually, one
considers the Maslov index as an invariant of a single loop in
$\L(n)$; this corresponds to taking $\lambda'$ to be a constant loop.
Let $L,L' \subset (M^{2n},\o)$ be two Lagrangian submanifolds. A loop
in $\PP$ is a map $u: S^1 \times I \longrightarrow M$ with boundary
values in $L$ resp. $L'$. After choosing a symplectic trivialization of
$u^*TM$ one obtains two loops $\lambda(s) = TL_{u(s,0)}$,
$\lambda'(s) = TL'_{u(s,1)}$ in $\L(n)$. The Maslov indices of such
loops determine a class $\chi \in H^1(\PP;\Z)$.

The {\em Maslov index for paths} assigns a half-integer
$\mu(\lambda,\lambda') \in \half\Z$ to any pair of paths
$\lambda,\lambda': [a;b] \longrightarrow \L(n)$. It is
a generalization of the ordinary Maslov index, to which it
reduces if both paths are closed, and has the following
basic properties:

\begin{romanlist}
\item \label{item:homotopy-invariance}
$\mu(\lambda,\lambda')$ depends on $\lambda,\lambda'$ only
up to homotopy with fixed endpoints.

\item \label{item:base-change} The Maslov index remains
the same if one conjugates both $\lambda$ and $\lambda'$ by
a path $\Psi: [a;b] \longrightarrow \Symp(2n,\R)$.

\item \label{item:additivity}
$\mu$ is additive under concatenation (of pairs of paths).

\item \label{item:reversal}
$\mu(\lambda,\lambda') = - \mu(\lambda',\lambda)$.

\item \label{item:clean}
$\mu(\lambda,\lambda')$ vanishes if the dimension of $\lambda(s) \cap
\lambda'(s)$ is constant.

\item \label{item:modone}
$\mu(\lambda,\lambda') \equiv \half \dim(\lambda(a) \cap \lambda'(a)) -
\half \dim(\lambda(b) \cap \lambda'(b)) \mymod 1.$
\end{romanlist}

Take a path $[a;b] \longrightarrow \PP$ from $\gamma_{x_-}$ to $\gamma_{x_+}$,
where $x_-,x_+ \in L \cap L'$. Such a path is given by a map
$u: [a;b] \times I \longrightarrow M$ with suitable boundary conditions. Let
$E = u^*(TM,\o)$. Choose a Lagrangian subbundle $F \subset E$ such
that $F_{(a,t)} = TL_{x_-}$ and $F_{(b,t)} = TL_{x_+}$ for all $t$, and $F_{(s,0)} =
TL_{u(s,0)}$ for all $s$. After choosing a trivialization of $E$, one
obtains two paths $\lambda,\lambda': [a;b] \longrightarrow \L(n)$,
namely $\lambda(s) = F_{(s,1)}$ and $\lambda'(s) = TL'_{u(s,1)}$.
Properties \ref{item:homotopy-invariance} and \ref{item:base-change}
ensure that
\[
I(u) \stackrel{\mathrm{def}}{=} \mu(\lambda,\lambda')
\]
is independent of the trivialization and of the choice of $F$. $I(u)$
depends on $u$ only up to homotopies which keep the endpoints
$\gamma_{x_{\pm}}$ fixed. It is also additive under concatenation.
Moreover, if $u$ and $u'$ are two paths with the same endpoints,
one has $I(u) - I(u') = \chi(v)$, where
$\chi \in H^1(\PP;\Z)$ is the class defined above and $v$ is the loop
in $\PP$ obtained by gluing $u$ and $u'$ at both endpoints. It
follows that if $\chi = 0$ one can find numbers $i(\gamma_x)
\in \half\Z$ for every $x \in L \cap L'$, such that
\[
I(u) = i(\gamma_{x_-}) - i(\gamma_{x_+})
\]
for every path $u$ from $\gamma_{x_-}$ to $\gamma_{x_+}$. Because
of property \ref{item:modone} one can also arrange that
\begin{equation} \label{eq:modone-one}
i(\gamma_x) \equiv \half \dim(TL_x \cap TL'_x) \mymod 1.
\end{equation}
Numbers $i(\gamma_x)$ with these two properties are called a
{\em coherent choice of indices} for $(L,L')$.

Take $H_-,H_+ \in \H$, and let $v: [a;b] \longrightarrow \PP$ be a
path from a point $\gamma_- \in Z(H_-)$ to a point $\gamma_+ \in
Z(H_+)$. A slight extension of the construction above associates to
such a path a number $I_{H_-,H_+}(v) \in \half\Z$. The details are
as follows: choose a trivialization of $E = v^*TM$ and a Lagrangian
subbundle $F \subset E$ such that
\[
F_{(a,t)} = \phi_t^{H_-}(TL_{\gamma_-(0)}), \quad
F_{(b,t)} = \phi_t^{H_+}(TL_{\gamma_+(0)}), \quad
F_{(s,0)} = TL_{v(s,0)}.
\]
Again, one obtains two paths in $\L(n)$: $\lambda(s) = F_{(s,1)}$
and $\lambda'(s) = TL'_{v(s,1)}$. $I_{H_-,H_+}(v)$ is defined as the
Maslov index of these paths.
If $\chi = 0$, one can find numbers $i_H(\gamma) \in \half\Z$ for any
$H \in \H$ and $\gamma \in Z(H)$, such that
\[
I_{H_-,H_+}(v) = i_{H_-}(\gamma_-) - i_{H_+}(\gamma_+)
\]
for all $v,H_-,H_+$, and
\begin{equation} \label{eq:modone-two}
i_{H}(\gamma) \equiv \half\dim(D\phi_1^H(TL_{\gamma(0)}) \cap
TL'_{\gamma(1)}) \mymod 1.
\end{equation}
Moreover, given a coherent choice of indices $i(\gamma_x)$, one
can choose the $i_H(\gamma)$ in such a way that
$i_H(\gamma_x) = i(\gamma_x)$ for $H = 0$.

Now let $\mo$ be the interior of a compact symplectic manifold
with contact type boundary. We assume that $L,L'$ are compact,
and that $[\alpha] \in H^1(\PP;\R)$ and $\chi \in H^1(\PP;\Z)$
vanish. Fix a coherent choice of indices for $(L,L')$, and extend
that choice to more general numbers $i_{H}(\gamma)$ as above. For
$H \in \H^{\reg}$ and $k \in \Z$, let $CF_k(H) \subset CF(H)$ be
the subgroup generated by those $\gamma \in Z(H)$ such that
$i_H(\gamma) = k$; it follows from \eqref{eq:modone-two} that
$i_H(\gamma)$ is always integral if $H \in \H^{\reg}$. Choose a
$\J \in \JJ^{\reg}(H)$. An index theorem due to Floer
\cite{floer88e} shows that $\partial(H,\J)$ has degree $-1$
with respect to the grading of $CF(H)$ which we have introduced.
Hence one obtains a grading of $HF(L,L',H,\J)$. This grading
is compatible with the canonical isomorphisms between these
groups for different $(H,\J)$.

One case when $\chi$ vanishes is when the first Chern class
of $\mo$ is zero and $H^1(L) = H^1(L') = 0$. This shows that
the grading of Floer homology exists in the situation described
in section \ref{sec:outline}.

Clearly, a choice of grading for $HF(L,L')$ also induces
a grading of all local Floer homology groups. As in the
previous section, assume that $N = L \cap L'$ has finitely
many path components $C_1, \dots, C_r$ which are open in
$N$. Then one obtains a filtration of the chain complex
$(CF_*(H),\partial(H,\J))$, for suitable $(H,\J)$, and the
homology of successive quotients is the local Floer homology
$HF_*^{\loc}(L,L';C_j)$. Therefore there is a spectral
sequence which converges to $HF_*(L,L')$, with
\begin{equation} \label{eq:general-ss}
E^1_{pq} = HF^{\loc}_{p+q}(L,L';C_p).
\end{equation}
Now assume that $L$ and $L'$ have clean intersection. It follows
from property \ref{item:clean} that for any coherent choice of
indices the function $x \longrightarrow i(\gamma_x)$ is locally
constant on $L \cap L'$. Let $i(C_j)$ be the value of this function
on $C_j$, and $i'(C_j) = i(C_j) - \half \dim C_j$ (equation
\eqref{eq:modone-one} implies that the $i'(C_j)$ are integral).
Theorem \ref{th:local-floer-homology} has the following graded
version:

\begin{theorem} \label{th:graded-pozniak}
$HF^{\loc}_*(L,L';C) \iso H_{*-i'(C)}(C;\Z/2)$.
\end{theorem}

Given this, one obtains the spectral sequence used in section
\ref{sec:outline} as a special case of \eqref{eq:general-ss}.
We will not prove Theorem \ref{th:graded-pozniak} but only the
case corresponding to Proposition \ref{th:small-pozniak}. To do this,
introduce local coordinates around $C$ as in the proof of that Proposition,
and take $H$ as in \eqref{eq:local-h}. If $h$ is sufficiently small,
the subset of $Z(H)$ which consists of paths near $C$ contains only the constant
paths $\gamma_{x_0}, \gamma_{x_1}$ at $x_i = (z_i,0,0,0)$, where $z_0$ and $z_1$
are the minimum and maximum of $h$. We must prove that
\begin{equation} \label{eq:change-of-index}
i_H(\gamma_{x_0}) = i'(C), \quad
i_H(\gamma_{x_1}) = i'(C)+1.
\end{equation}
By definition $i_H(\gamma_{x_0})$ has the following
property: take a map $u: I^2 \longrightarrow M$ such that
$u(0,t) = u(1,t) = x_0$, $u(s,0) \in L$ and $u(s,1) \in L'$
for all $s,t$. Then
\[
i(\gamma_{x_0}) - i_H(\gamma_{x_0}) = I_{0,H}(u).
\]
$u$ can be chosen to be the constant map at $x_0$. The local
coordinates which we are using provide a trivialization
of $u^*TM$. To compute $I_{0,H}(u)$ one has to choose a subbundle
$F \subset u^*TM$ with certain properties: one possible
choice is
\[
F_{(s,t)} = D\phi^H_{st}(TL_{x_0})
= \{ (r,r\cdot s\cdot t\cdot h''(z_0)) : r \in \R \} \times \R
\times 0.
\]
$I_{0,H}(u)$ is defined as the Maslov index of the paths
$\lambda(s) = F_{(s,1)}$, $\lambda'(s) = TL_{x_0}' =
\R \times 0 \times 0 \times \R$. Using the definition in
\cite{robbin-salamon93} and the fact that $h''(z_0)>0$,
one obtains $\mu(\lambda,\lambda') = \half$. Therefore
$i_H(\gamma_{x_0}) = i(\gamma_{x_0}) - \half = i'(C)$.
The same argument can be used to prove the second part of
\eqref{eq:change-of-index}.

\section{Geodesics \label{sec:tangent}}

This section summarizes the classical relationship between geodesics and
Lagrangian intersections. Let $(P,g)$ be a compact Riemannian manifold
and $c: I \longrightarrow P$ a geodesic. For $r \in I$, let $m(c,r) \in \Z$
be the multiplicity of $c(0)$ and $c(r)$ as conjugate points along $c$.
The energy and Morse index of $c$ are defined by
\begin{align}
\notag e(c) &= \half g(\dot{c}(0),\dot{c}(0)),\\
\label{eq:morse} m(c) &= \sum_{0<r<1} m(c,r) + \half m(c,1).
\end{align}
$m(c)$ is not necessarily an integer; we have adjusted the contribution of the
endpoints to suit Robbin-Salamon's conventions for the Maslov index.

Let $\mo$ be the tangent bundle $TP$ together with the symplectic form
obtained by identifying $TP \iso \cotangent P$. Let $(\phi_t)_{t \in \R}$
be the geodesic flow on $M$, that is, the Hamiltonian flow of $H(\xi) =
\half g(\xi,\xi)$. Choose two points $p,p' \in P$, and consider the
Lagrangian submanifolds $L = \phi_{-1}(TP_p), L' = TP_{p'} \subset M$.
Their intersection points correspond to geodesics from $p'$ to $p$. More
precisely, a point $\xi \in L'$ lies in $N = L \cap L'$ iff the unique geodesic
$c_{\xi}: I \longrightarrow P$ with $\dot{c}_{\xi}(0) = \xi$ satisfies
$c_{\xi}(1) = p$. The numbers $m(c_{\xi},r)$ can be written in terms
of the derivative of $\phi$:
\begin{equation} \label{eq:conjugate}
m(c_{\xi},r) = \dim\left(\Lambda_{\xi} \cap [D\phi_{-r}(\Lambda)]_{\xi}\right),
\end{equation}
where $\Lambda \subset TM$ is the vertical part of $TM$, that is, the tangent bundle
along the fibres of the projection $M \longrightarrow P$. In particular
$m(c_{\xi},1) = \dim (TL_{\xi} \cap TL'_{\xi})$.
Therefore $L$ and $L'$ have clean intersection iff $N$ is a
submanifold and $\dim N = m(c_{\xi},1)$ for all $\xi \in N$; the last
condition means that every Jacobi field along $c_{\xi}$ with vanishing
boundary values comes from a geodesic variation of $c_{\xi}$ which
leaves the endpoints fixed.

It is easy to see in the present case both $[\alpha] \in H^1(\PP;\R)$
and the class $\chi \in H^1(\PP;\Z)$ defined in the previous section vanish.
Hence one can choose an action functional $a: \PP \longrightarrow \R$, and
a coherent choice of indices $i(\gamma_{\xi})$. Both are not unique; the
following Proposition holds for one particular choice.

\begin{prop} \label{th:geodesics}
Let $\gamma_{\xi} \in \PP$ be the constant path at
a point $\xi \in L \cap L'$, and $c_{\xi}$ the corresponding geodesic.
Then $a(\gamma_{\xi}) = e(c_{\xi})$ and $i(\gamma_{\xi}) = m(c_{\xi})$.
\end{prop}

We begin by considering a slightly more general situation.

\begin{lemma} \label{th:computational}
Let $\mo$ be a symplectic manifold, and $L_0,L'_0 \subset M$
a pair of Lagrangian submanifolds such that
\begin{normallist}
\item there is a $\theta \in \Omega^1(M)$ with $d\theta = \o$ and
$\theta|L_0 = 0$, $\theta|L'_0 = 0$;
\item there is a Lagrangian subbundle $\Lambda \subset TM$ with
$\Lambda|L_0 = TL_0$, $\Lambda|L'_0 = TL'_0$.
\end{normallist}
Take a proper function $H \in \smooth(M;\R)$ with Hamiltonian vector
field $X$, and let $(\phi_t)$ be its flow. Set $L = \phi_{-1}(L_0)$,
$L' = L_0'$. Then the two classes $[\alpha] \in H^1(\P(L,L');\R)$ and
$\chi \in H^1(\P(L,L');\Z)$ vanish, and for a suitable choice of action
$a$ and indices $i$, the following holds:
\begin{primelist}
\item $a(\gamma_x) = -H(x) + \int_I (i_X\theta)(\phi_t(x))\, dt$
for all $x \in L \cap L'$.
\item \label{item:explicit-maslov} Take $x \in L \cap L'$ and choose a
symplectic isomorphism $TM_x \iso \R^{2n}$. Then
\begin{equation} \label{eq:explicit-maslov}
i(\gamma_x) = \mu(\lambda_x,\lambda'_x) - \half \dim L,
\end{equation}
where $\lambda_x,\lambda'_x: I \longrightarrow \L(n)$ are given by
$\lambda_x(r) = \Lambda_x$ and $\lambda_x'(r) = [D\phi_{-r}(\Lambda)]_x$.
\end{primelist}
\end{lemma}

\proof We prove only the statement \ref{item:explicit-maslov}
and leave the rest to the reader. Take two points $x_-,x_+ \in
L \cap L'$ and a map $u: [a;b] \times I \longrightarrow M$ which
corresponds to a path from $\gamma_{x_-}$ to $\gamma_{x_+}$ in
$\PP$. In order to compute $I(u)$ one has to choose a trivialization
of $E = u^*TM$ and a Lagrangian subbundle $F \subset E$ with
certain properties. One suitable choice is
$F_{(s,t)} = [D\phi_{-1}(\Lambda)]_{u(s,t)}$. $I(u)$ is the
Maslov index of the pair $(\lambda,\lambda')$ given by
\[
\lambda(s) = F_{(s,1)} = [D\phi_{-1}(\Lambda)]_{u(s,1)}, \quad
\lambda'(s) = TL'_{(s,1)} = \Lambda_{u(s,1)}.
\]
Consider another Lagrangian subbundle $F' \subset E$ defined by
\[
F'_{(s,t)} = [D\phi_{s-1}(\Lambda)]_{u(s,t)}.
\]
For any path $\alpha$ in $[a;b] \times I$, we denote by $\tilde{\mu}(\alpha)$ the
Maslov index of $(F,F')$ along $\alpha$, that is, the Maslov
index of $r \longmapsto (F_{\alpha(r)},F'_{\alpha(r)})$.
For instance, the expression for $I(u)$ given above says that
$I(u) = \tilde{\mu}(\alpha_2)$ where
$\alpha_2: [a;b] \longrightarrow [a;b] \times I$ is the path $\alpha_2(r) = (r,1)$.
Now take the other three sides of the boundary of $[a;b] \times I$:
$\alpha_1(r) = (r,0)$ $(a \leq r \leq b)$ and $\alpha_3(r) = (a,r)$,
$\alpha_4(r) = (b,r)$ $(0 \leq r \leq 1)$. Because of the additivity
and homotopy invariance of the Maslov index for paths,
\begin{equation} \label{eq:additivity}
I(u) = \tilde{\mu}(\alpha_2) = -\tilde{\mu}(\alpha_3) + \tilde{\mu}(\alpha_1) + \tilde{\mu}(\alpha_4).
\end{equation}
Since $F$ and $F'$ agree over $[a;b] \times 0$, $\tilde{\mu}(\alpha_1) = 0$.
$\tilde{\mu}(\alpha_3)$ and $\tilde{\mu}(\alpha_4)$ are independent of $u$; they
depend only on $x_-$ and $x_+$, respectively. After changing the
trivialization of $E$ by $D\phi_{1-s}$ one sees that $\tilde{\mu}(\alpha_3) =
\mu(\lambda_{x_-}',\lambda_{x_-})$, and therefore (by property \ref{item:reversal}
of $\mu$) $-\tilde{\mu}(\alpha_3) = \mu(\lambda_{x_-},\lambda_{x_-}')$. Similarly,
$-\tilde{\mu}(\alpha_4) = \mu(\lambda_{x_+},\lambda_{x_+}')$. Equation \eqref{eq:additivity}
says that \eqref{eq:explicit-maslov} is a coherent choice of indices.
The constant $\half \dim L$ has been subtracted in order to
fulfil the integrality criterion \eqref{eq:modone-one}. \qed

\proof[Proof of Proposition \ref{th:geodesics}] Let $\theta \in \Omega^1(M)$
be the form corresponding to the canonical one-form on $\cotangent P$,
$\Lambda \subset TM$ the vertical subbundle, and $H(\xi) = \half g(\xi,\xi)$.
$L_0 = TP_p$ and $L_0' = TP_{p'}$ satisfy the conditions of Lemma \ref{th:computational}.
Using the first part of that Lemma and the fact that $i_X\theta = 2H$, one obtains
\[
a(\gamma_{\xi}) = -H(\xi) + \int_I (i_X\theta)(\phi_t(x)) \; dt = H(\xi)
\]
for any $\xi \in L \cap L'$.
Choose a symplectic isomorphism $TM_{\xi} \iso \R^{2n}$ induced by an
isomorphism $TP_{p'} \iso \R^n$ and by the Levi-Civita connection.
Then the paths $\lambda_{\xi},\lambda_{\xi}'$ defined in Lemma \ref{th:computational}
are of the following form: $\lambda_{\xi}(r) = \R^n \times 0$, and
$\lambda_{\xi}'(r) = A(r)^{-1}(\R^n \times 0)$, where
$A: [0;1] \longrightarrow \Symp(2n,\R)$ satisfies a differential
equation
\begin{equation} \label{eq:ode}
\dot{A}(r) = \begin{pmatrix} 0 & R(r) \\ \One & 0 \end{pmatrix} A(r), \quad
A(0) = \One
\end{equation}
for some family $R(r)$ of symmetric $n\times n$-matrices obtained from the
curvature tensor of $(P,g)$. This is just the equation for Jacobi fields,
written as a first order equation. In view of \eqref{eq:morse} and \eqref{eq:conjugate},
the proof of Proposition \ref{th:geodesics} is completed by applying the
following property of the Maslov index for paths:

\begin{lemma} Let $R(r)$, $0 \leq r \leq 1$, be a family of symmetric $n \times n$
matrices, and let $A(r)$ be the solution of \eqref{eq:ode}. Consider paths
$\lambda,\lambda': [0;1] \longrightarrow \L(n)$ given by $\lambda(r) = \R^n \times 0$,
$\lambda'(r) = A(r)^{-1}(\R^n \times 0)$. Their Maslov index is
\begin{multline*}
\mu(\lambda,\lambda') = \half \dim(\lambda(0) \cap \lambda'(0)) +
\sum_{0<r<1} \dim(\lambda(r) \cap \lambda'(r)) + \\ + \half
\dim(\lambda(1) \cap \lambda'(1)).
\end{multline*}
\end{lemma}

This property can be deduced easily from the definition of $\mu$ given
in \cite{robbin-salamon93}.

\section{Generalized Dehn twists \label{sec:dehn-twists}}

This section contains the definition of the maps $\tau_L$. The following
elementary fact will be used several times:

\begin{lemma} \label{th:elementary}
Let $\mo$ be a symplectic manifold, $H \in \smooth(M,\R)$ and
$\Psi \in \smooth(\R,\R)$. The Hamiltonian flows of $H$ and
$\Psi(H)$ are related by
\[
\phi^{\Psi(H)}_t(x) = \phi^H_{t\Psi'(H(x))}(x). \qed
\]
\end{lemma}

Let $\eta$ be the standard symplectic form on $\cotangent S^2$, and
$S^2 \subset \cotangent S^2$ the zero-section. Its complement
$\cotangent S^2\setminus S^2$ carries a Hamiltonian circle action
$\sigma$ with moment map $\mu(\xi) = |\xi|$ (the length function
with respect to the standard metric). To see that this is a
circle action, recall that if we identify $\cotangent S^2 = TS^2$ then
the flow induced by $\half\mu^2$ is the geodesic flow. By Lemma
\ref{th:elementary}, $\mu$ itself induces the {\em normalized
geodesic flow} which transports any nonzero tangent vector $\xi$
with unit speed along the geodesic emanating from it, irrespective
of what $|\xi|$ is. Since all geodesics of length $2\pi$ are closed,
this is a circle action. $\sigma$ does not extend continuously
over the zero-section, with one exception: since any geodesic of
length $\pi$ on $S^2$ connects two opposite points, $\sigma(-1)$
is the restriction of the antipodal involution $A: \cotangent S^2
\longrightarrow \cotangent S^2$ to $\cotangent S^2 \setminus S^2$.

Take a function $\psi \in \smooth(\R,\R)$ such that $\psi(t) +
\psi(-t) = 2\pi$ for all $t$, and $\psi(t) = 0$ for $t \gg 0$.
Let $\tau: \cotangent S^2 \longrightarrow \cotangent S^2$ be
the map defined by
\[
\tau(\xi) = \begin{cases}
\sigma(e^{i\psi(|\xi|)})(\xi) & \xi \notin S^2,\\
A(\xi) & \xi \in S^2.\\
\end{cases}
\]
$\tau$ is smooth and symplectic, and it is the identity outside a
compact subset. The third property is obvious; to prove the first
two, consider
\[
(A \circ \tau)(\xi) = \sigma(e^{i(\psi(|\xi|)-\pi)})(\xi).
\]
Take a function $\Psi \in \smooth_c(\R,\R)$ with $\Psi'(t) =
\psi(t) - \pi$. Since $\psi-\pi$ is odd, $\Psi$ is even, and
hence $\xi \longmapsto \Psi(|\xi|)$ is smooth on all of $\cotangent S^2$.
Lemma \ref{th:elementary} shows that $A \circ \tau$ is the
time-one map of the Hamiltonian flow of $\Psi(|\xi|)$. In particular, it is smooth and
symplectic, and therefore so is $\tau$. We call $\tau$ a {\em model generalized Dehn twist}.

Let $\mo$ be a symplectic four-manifold containing a Lagrangian
two-sphere $L$. By a theorem of Weinstein, there is a symplectic
embedding $f: \disc \longrightarrow M$ of the disc bundle
$\disc = \{ \xi \in \cotangent S^2 \suchthat |\xi| < \epsilon\}$
into $M$, for some $\epsilon>0$, such that $f(S^2) = L$. Let
$\tau$ be the model generalized Dehn twist associated to a function
$\psi$ such that $\psi(t) = 0$ for all $t>\epsilon/2$. Then one can
define a symplectic automorphism $\tau_L$ of $M$ by
\[
\tau_L(x) = \begin{cases}
f \tau f^{-1}(x) & x \in \im(f),\\
x & \text{otherwise}.
\end{cases}
\]
We call such a map $\tau_L$ a generalized Dehn twist along $L$.

\begin{lemma} \label{th:dehn-twist-well-defined}
The symplectic isotopy class of $\tau_L$ is independent of the
choice of $f$ and $\psi$. Moreover, if $L$ and $L'$ are
Lagrangian isotopic then $\tau_L$ and $\tau_{L'}$ are symplectically
isotopic. \end{lemma}

The independence of $\psi$ can be proved by an explicit isotopy.
Next, consider two embeddings $f,f': \disc \longrightarrow M$
with $f(S^2) = f'(S^2) = L$. If $f$ can be deformed to
$f'$ through symplectic embeddings which map $S^2$ to $L$
then the corresponding generalized Dehn twists are
symplectically isotopic. The same holds if $f$ can be
deformed to $f'$ after making $\epsilon$ smaller. Such
a deformation of the germs of $f,f'$ exists iff the
restrictions $f|S^2,f'|S^2: S^2 \longrightarrow L$ are
differentiably isotopic. Since $\Diff^+(S^2)$ is path-connected,
this holds iff $f$ and $f'$ induce the same orientation
of $L$. To complete the proof that $\tau_L$ is independent of
the choice of embedding, it is enough to find two examples $f,
f'$ which induce opposite orientations of $L$ but define the
same generalized Dehn twist, and that is easy: take an arbitrary
$f$ and set $f' = f \circ A$. Finally, it is clear that the
symplectic isotopy class of $\tau_L$ depends on $L$ only up to
ambient symplectic isotopy. However, that is the same as Lagrangian
isotopy.

An inspection of the proof which we have just given shows that
$\tau_L$ is well-defined up to Hamiltonian isotopy. We do not need
this sharper statement here.

\begin{lemma} \label{th:order-two}
Let $\tau_L$ be a generalized Dehn twist along a Lagrangian two-sphere
$L$. Then the square $\tau_L^2$ is differentiably isotopic to the identity.
\end{lemma}

\proof We use the model
$\cotangent S^2 = \{ (u,v) \in \R^3 \times \R^3 \suchthat
|u|=1 \text{ and } \leftsc u,v \rightsc = 0 \}$,
in which $\eta = \sum_i dv_i \wedge du_i$. For $x \in \R^3 \setminus 0$ and
$t \in \R$, let $R^t(x) \in SO(3)$ be the rotation with axis $x/|x|$ and
angle $t$. Then
\[
\sigma(e^{it})(u,v) = (R^t(u \times v)u,R^t(u \times v)v).
\]
Consider the following one-parameter family $\sigma^{(s)}$, $0 \leq s \leq 1$,
of smooth circle actions on $\cotangent S^2 \setminus S^2$:
\[
\sigma^{(s)}(e^{it})(u,v) = (R^t(su + (1-s)u \times v)u, R^t(su + (1-s)u \times v)v).
\]
$\sigma^{(0)} = \sigma$. On the other hand, $\sigma^{(1)}$ is the action of $S^1$ by rotation in each
fibre of $T^*S^2$ and extends smoothly to the zero-section $S^2$. The square of a
model generalized Dehn twist is
\[
\tau^2(\xi) = \begin{cases}
\sigma(e^{2i\psi(|\xi|)})(\xi) & \xi \notin S^2,\\
\xi & \xi \in S^2.
\end{cases}
\]
We can assume that $\psi(t) = \pi$ for small $|t|$; then $\tau^2$ is the identity in
a neighbourhood of the zero-section.
Replacing $\sigma$ by $\sigma^{(s)}$ defines a differentiable isotopy from $\tau^2$ to
$T(\xi) = \sigma^{(1)}(e^{2i\psi(|\xi|)})(\xi)$, and this can be deformed to the
identity by changing $\psi$ to $s\psi$ for $0 \leq s \leq 1$. This isotopy from
$\tau^2$ to the identity is local in the sense that if $\tau = \id$ outside
$\disc$ for some $\epsilon>0$ then the same holds for the isotopy.
This implies the Lemma as stated. \qed

\begin{remarks} \label{rem}
\begin{remarkslist}
\item \label{item:naturality}
Let $L \subset M$ be a Lagrangian two-sphere and $f$ a symplectic
automorphism of $M$. It follows immediately from the definition of
generalized Dehn twists that $\tau_{f(L)} = f \tau_L f^{-1}$.

\item The definition of a generalized Dehn twist is sensitive to
the sign of $\o$. A generalized Dehn twist along $L$ as a submanifold of $(M,-\o)$
is the inverse of a generalized Dehn twist along $L \subset (M,\o)$.

\item Let $\Aut^c(\cotangent S^2,\eta)$ be the group of those symplectic
automorphisms of $\cotangent S^2$ which are equal to the identity outside
a compact subset, and $[\tau] \in \pi_0(\Aut^c(\cotangent S^2,\eta))$ the
class containing all model generalized Dehn twists. Corollary \ref{th:corollary}
implies that $[\tau]$ has infinite order. It can be shown \cite{seidel98c}
that $[\tau]$ generates $\pi_0(\Aut^c(\cotangent S^2,\eta))$, and that the
higher homotopy groups are trivial.

\item The definition of a model generalized Dehn twists extends in a
straightforward way to the cotangent bundle of $S^n$ for all $n$. Using this as a local
model, one can define generalized Dehn twists associated to Lagrangian embeddings of $S^n$
into $2n$-dimensional symplectic manifolds. For $n = 1$ these are just the ordinary
positive Dehn twists along a curve on a surface.
\end{remarkslist}
\end{remarks}
\section{Proof of Theorem \ref{th:main}}

Let $\mo$ and $L_1,L_2,L_3$ be as in that Theorem. Fix some $r \in \N$. One can find
a symplectic embedding $f: \disc \longrightarrow M$ for some $\epsilon>0$, such
that $f(S^2) = L_2$, $f^{-1}(L_1) = \cotangent_x S^2 \cap \disc$ and
$f^{-1}(L_3) = \cotangent_{A(x)} S^2 \cap \disc$ for some $x \in S^2$. After
rescaling $\o$ if necessary, one can assume that $\epsilon = 2\pi r$. Let
$\tau_{L_2}$ be the generalized Dehn twist along $L_2$ defined using the
embedding $f$ and some function $\psi$, and $L_1^{(r)} = \tau_{L_2}^{2r}(L_1)$.
$L_1^{(r)} \cap L_3$ is contained in $\im(f)$, and
$
f^{-1}(L_1^{(r)} \cap L_3) =
\tau^{2r}(\cotangent_x S^2) \cap \cotangent_{A(x)}S^2
= \{ \xi \in \cotangent_{A(x)}S^2 \suchthat 2r\psi(|\xi|) + \pi \in 2\pi\Z\}$,
where $\tau$ is the model generalized Dehn twist determined by $\psi$.
Now assume that $\psi$ satisfies
\[
\begin{cases}
\psi'(t) \leq 0 & \text{for all } t,\\
\psi(t) = \pi - t/2r & \text{for } 0 \leq t \leq \delta = 2\pi(r - \quarter),\\
\psi(t) = 0 & \text{for } t \geq 2\pi r
\end{cases}
\]
Then $L_1^{(r)} \cap L_3$ is the disjoint union of $r$ circles $C_1, \dots, C_r$, where
\[
f^{-1}(C_j) = \{\xi \in \cotangent_{A(x)}S^2 \suchthat \psi(|\xi|) = \textstyle{\frac{2j-1}{2r}}\pi\}.
\]
Note that all $C_j$ lie in $f(U)$, where $U = D_{\delta}(\cotangent S^2)
\subset D_\epsilon(\cotangent S^2)$. This is important because
$\tau^{2r}(\xi) = \sigma(e^{-i|\xi|})(\xi)$ for all $\xi \in U$.
Since $\sigma$ is defined by normalizing the geodesic flow $(\phi_t)$, this means
that $\tau^{2r}|U = \phi_{-1}|U$. Setting
$L = \phi_{-1}(\cotangent_xS^2)$ and $L' = \cotangent_{A(x)}S^2$, we have shown that
\[
f^{-1}(L_1^{(r)}) \cap U = L \cap U, \quad f^{-1}(L_3) \cap U = L' \cap U.
\]
This makes it possible to apply the results of section \ref{sec:tangent}.
First of all, the intersection points of $L_1^{(r)}$ and $L_3$ are the intersection
points of $L$ and $L'$ inside $U$, and these correspond to geodesics from $A(x)$
to $x$ of length $\leq \delta$. More precisely, the circle $C_j$
corresponds to the one-parameter family of geodesics which wind $j - \half$ times
around $S^2$. In section \ref{sec:tangent} we have given a criterion, in terms of
Jacobi fields, for $L$ and $L'$ to have clean intersection. This is satisfied is
the present case. Hence $L_1^{(r)}$ and $L_3$ also have clean intersection.
To compute the relative action and index of two intersection points $x_-,x_+ \in
L_1^{(r)} \cap L_3$ one can use a path in $\P(L_1^{(r)},L_3)$ whose image lies
inside $f(U)$. Therefore the relative action and index coincide with those of $f^{-1}(x_-),
f^{-1}(x_+)$ as intersection points of $L$ and $L'$. Using Proposition \ref{th:geodesics},
one obtains that the action $a_j \in \R$ of a constant path at a point of $C_j$ satisfies
\[
a_j - a_{j-1} = \textstyle
\frac{\pi^2}{2}((2j-1)^2 - (2j-3)^2) > 0.
\]
The Morse index of a geodesic from $A(x)$ to $x$ which winds $j - \half$ times
around $S^2$ is $2j - \frac{3}{2}$ (it has $2j$ conjugate points on it, including
both endpoints, and all of them have multiplicity one; according to our definition,
one endpoint does not contribute at all, while the other contributes $\half$). Therefore
\[
i'(C_j) - i'(C_{j-1}) = 2.
\]
This completes the computations necessary to apply Pozniak's spectral
sequence, as described in section \ref{sec:outline}. Note that since $L_1^{(r)}$
and $L_3$ are orientable and intersect in a union of circles, we have provided
proofs of the basic results underlying the spectral sequence
(see Proposition \ref{th:small-pozniak} and the discussion following
Theorem \ref{th:graded-pozniak}).

\section{A family of examples \label{sec:examples}}

A {\em configuration of Lagrangian two-spheres} in a symplectic four-manifold
is a finite collection of Lagrangian two-spheres any two of which
intersect transversely. An $(A_m)$-configuration, for $m \geq 1$, consists
of $m$ Lagrangian two-spheres $L_1,\dots,L_m$ such that
\begin{equation} \label{eq:am}
|L_i \cap L_j| = \begin{cases} 1 & i-j = \pm 1,\\ 0 & |i-j| \geq 2. \end{cases}
\end{equation}

\begin{prop} \label{th:am}
Let $(H,\o)$ be the affine hypersurface $z_1^2 + z_2^2 = z_3^{m+1} + \half$ in $\C^3$, equipped
with the standard symplectic form. For any $m$, $(H,\o)$ contains an $(A_m)$-configuration
of Lagrangian two-spheres.
\end{prop}

\proof The projection $\pi: H \longrightarrow \C^2$ onto $(z_1,z_2)$ is an
$(m+1)$-fold covering branched along $C = \{z_1^2 + z_2^2 = \half\} \subset \C^2$.
The covering group is generated by $\sigma(z_1,z_2,z_3) = (z_1,z_2,e^{2\pi i/(m+1)}z_3)$.
Let $\o_0$ be the standard symplectic form on $\C^2$.

\begin{lemma} \label{th:branched-cover}
Let $K \subset \C^2 \setminus C$ be a compact subset. There is a symplectic form $\o'$ on $H$
which is diffeomorphic to $\o$ and such that
\begin{equation} \label{eq:pullback}
\o' | \pi^{-1}(U) = \pi^*(\o_0|U)
\end{equation}
for some neighbourhood $U \subset \C^2$ of $K$.
\end{lemma}

\proof[Proof of Lemma \ref{th:branched-cover}] $K$ is contained in the open
subset $U = \{ z \in \C^2 \suchthat \epsilon < |z_1^2 + z_2^2 - \half|^{1/(m+1)}
< \epsilon^{-1} \}$ for sufficiently small $\epsilon>0$. Choose
a function $\beta \in \smooth(\Rgeq,\R)$ such that $\beta(r) \leq 1$ for all
$r$, $\beta(r) = 0$ for $r \leq \epsilon/2$ or $r \geq 2\epsilon^{-1}$,
$\beta(r) = 1$ for $\epsilon \leq r \leq \epsilon^{-1}$, and
$\int_0^\infty r \beta(r)\; dr = 0$. Set
\[
\o' = \o - \beta(|z_3|) ({\textstyle\frac{i}{2}} dz_3 \wedge d\bar{z}_3).
\]
By definition $\o'$ satisfies \eqref{eq:pullback}. Moreover, it is a
symplectic form which is compatible with the complex structure; it agrees
with $\o$ outside a compact subset; and the difference $\o'-\o$
represents the trivial class in $H^2_c(H;\R)$. By a familiar argument
it follows that $\o$ and $\o'$ are diffeomorphic. \qed

Now consider the two-dimensional figure-eight map
\[
\R^3 \supset S^2 \stackrel{f}{\longrightarrow} \C^2 \setminus C,
\quad f(t_1,t_2,t_3) = (t_2(1+it_1),t_3(1+it_1)).
\]
$f$ is an immersion with one double point $0 = f(\pm 1,0,0)$ at which
the two branches meet transversely. Moreover, if $\gamma: [0;1] \longrightarrow
S^2$ is any path from $(1,0,0)$ to $(-1,0,0)$, $f(\gamma)$ is a loop in
$\C^2 \setminus C$ whose linking number with $C$ equals one. Let
$\tilde{f}: S^2 \longrightarrow H$ be a lift of $f$ to $H$; such a lift
exists because $f$ avoids the branch locus of $\pi$. The fact that
$f(\gamma)$ has linking number $1$ with $C$ implies that
\begin{equation} \label{eq:intersection}
\tilde{f}(-1,0,0) = \sigma(\tilde{f}(1,0,0)).
\end{equation}
Therefore $\tilde{f}$ is an embedding. Now consider the shifted embedding
$\sigma \circ \tilde{f}$. If $m \geq 2$, the images of $\tilde{f}$ and of
$\sigma \circ \tilde{f}$ do not have any intersection points except for
\eqref{eq:intersection}. The intersection at that point is modelled on the
self-intersection of $f$; hence it is transverse. A repetition of the
same argument shows that $L_1 = \im(\tilde{f}), L_2 = \sigma(L_1),
\dots, L_m = \sigma^{m-1}(L_1)$ is a family of smoothly embedded
two-spheres which intersect according to \eqref{eq:am}. Take a symplectic form
$\o'$ as in Lemma \ref{th:branched-cover} with $K = \im(f)$. Since $f$ is a
Lagrangian immersion with respect to $\o_0$, the submanifolds $L_1, \dots,
L_m$ are $\o'$-Lagrangian. This proves that $(H,\o')$ contains an $(A_m)$-configuration.
Since $\o'$ is diffeomorphic to $\o$, it follows that $(H,\o)$ contains one as well. \qed

\begin{remark}
\begin{remarkslist}
\item By leaving out some components, one sees that $(H,\o)$ contains an $(A_3)$-configuration
whenever $m \geq 3$. This configuration lies in the bounded subset \eqref{eq:affine}
if $R$ is sufficiently large. Taking $R$ large also ensures that the closure of \eqref{eq:affine} is
a symplectic manifold with contact type boundary. Both the symplectic class and the
first Chern class of \eqref{eq:affine} vanish, because it is an open subset of an affine hypersurface.
Hence \eqref{eq:affine} satisfies the conditions of Theorem \ref{th:main}.

\item The existence of $m$ smooth embedded two-spheres in $H$ satisfying
\eqref{eq:am} is a consequence of Brieskorn's resolution \cite{brieskorn66}.
The only new aspect of Proposition \ref{th:am} is that one can choose these
spheres to be Lagrangian.

\item A straightforward generalization of the proof given above produces an
$(A_m)$-configuration of Lagrangian $n$-spheres in the hypersurface
$z_1^2 + z_2^2 + \cdots + z_n^2 = z_{n+1}^{m+1} + \half$ for any $m,n$.
\end{remarkslist}
\end{remark}

\appendix
\section{Lagrangian surgery \label{sec:braiding}}

The aim of this Appendix is to relate generalized Dehn twists
to the {\em Lagrangian surgery} construction which has been
studied by Polterovich \cite{polterovich91} and others. As a
by-product we obtain the following result:

\begin{proposition} \label{th:braiding}
Let $L_1$ and $L_2$ be two Lagrangian
two-spheres in a symplectic four-manifold $\mo$. Assume
that they intersect transversely in a single
point. Then $\tau_{L_1} \tau_{L_2} \tau_{L_1}$ and
$\tau_{L_2} \tau_{L_1} \tau_{L_2}$ are symplectically
isotopic automorphisms.
\end{proposition}

In particular, an $(A_m)$-configuration in a symplectic
four-manifold defines a homomorphism from the braid group
$B_{m+1}$ to the group of symplectic isotopy classes of
automorphisms of the manifold. This holds e.g. for the
manifolds \eqref{eq:affine}.

We begin by recalling the definition of Lagrangian surgery.
Our exposition follows \cite{polterovich91} with some
modifications. Let $C \subset \R^2$ be a smooth embedded
curve with the following properties: $C$ is diffeomorphic
to $\R$; it coincides with
$(\R^+ \times 0) \cup (0 \times \R^-)$ outside a compact subset;
and there is no $x \in \R^2$ such that both $x$ and $-x$ lie in $C$.
Consider
\[
H = \{ (y_1 \cos t, y_1 \sin t, y_2 \cos t, y_2 \sin t)
\suchthat (y_1,y_2) \in C, \; t \in S^1 \} \subset \R^4.
\]
$H$ is an embedded surface diffeomorphic to $\R \times S^1$;
it is Lagrangian with respect to $\o = dx_1 \wedge dx_3 + dx_2 \wedge dx_4$;
and it coincides with $(\R^2 \times 0) \cup (0 \times \R^2)$
outside a compact subset. By choosing $C$ suitably, one can
arrange that the last-mentioned property holds outside an
arbitrarily small neighbourhood of $0 \in \R^4$. $H$ is
called a {\em Lagrangian handle}.

Now let $\mo$ be a symplectic four-manifold and $L_1,L_2 \subset M$ two compact Lagrangian
surfaces which intersect transversely and in a single point $x$. Choose a neighbourhood $U
\subset \R^4$ of $0$ and a Darboux chart $f: U \longrightarrow M$ such that $f(0) = x$,
$f^{-1}(L_1) = (\R^2 \times 0) \cap U$ and $f^{-1}(L_2) = (0 \times \R^2) \cap U$. Let $H$ be
a Lagrangian handle which agrees with $(\R^2 \times 0) \cup (0 \times \R^2)$ outside $U$.
Define a new Lagrangian submanifold $L \subset M$ by $L \cap f(U) = f(H)$ and $L \setminus
f(U) = (L_1 \cup L_2) \setminus f(U)$. $L$ is diffeomorphic to the connected sum of $L_1$ and
$L_2$. It is called the Lagrangian surgery of $L_1$ and $L_2$; we denote it by $L_1 \surgery
L_2$. One can show that this surgery is independent of all choices up to Lagrangian isotopy.

\begin{proposition} \label{th:surgery-one}
Assume that $L_2$ is a Lagrangian sphere.
Then $L_1 \surgery L_2$ is Lagrangian isotopic to
$\tau_{L_2}^{-1}(L_1)$.
\end{proposition}

\proof Let $\tau$ be the model generalized Dehn twist on
$\cotangent S^2$ defined using a function $\psi$ such that
$\psi' \leq 0$ everywhere, and
\[
\psi(t) = \pi-t \text{ for } t \leq \epsilon, \quad
\psi(t) > 0 \text{ for } \epsilon < t < 2\epsilon, \quad
\psi(t) = 0 \text{ for } t \geq 2\epsilon
\]
for some $\epsilon>0$. Choose a point $x \in S^2$, and set $L = \tau^{-1}(\cotangent_xS^2)$.
The exponential maps at $x$ and $A(x)$ induce symplectic isomorphisms
\begin{align*}
f_x: \cotangent B_\pi &\longrightarrow \cotangent S^2 \setminus \cotangent_{A(x)}S^2,\\
f_{A(x)}: \cotangent B_\pi &\longrightarrow \cotangent S^2 \setminus \cotangent_{x}S^2.
\end{align*}
Here $B_{\pi} \subset \R^2$ is the open ball of radius $\pi$. We will identify
$\cotangent B_{\pi}$ with $\R^2 \times B_\pi$. In these coordinates
\begin{align*}
f_x^{-1}(L) &=
\textstyle{\{ (p,-\frac{\psi(|p|)}{|p|}p) \suchthat p \in \R^2 \setminus 0 \},}\\
f_{A(x)}^{-1}(L) &=
\textstyle{\{ (p,-\frac{\pi-\psi(|p|)}{|p|}p) \suchthat p \in B_{2\epsilon}\}.}
\end{align*}
Let $V = B_{\epsilon} \times B_{\epsilon} \subset \cotangent B_{\pi}$. Then
\[
f_{A(x)}^{-1}(L) \cap V = \{(p,-p) \suchthat p \in B_{\epsilon}\}.
\]
It follows that $L$ can be deformed (by a symplectic isotopy which is trivial outside $f_{A(x)}(V)$)
into the Lagrangian submanifold $L' \subset \cotangent S^2$ defined by
\begin{align*}
f_x^{-1}(L') &=
\textstyle{\{ (\rho(\pi-\psi(|p|))p,-\frac{\psi(|p|)}{|p|}p) \suchthat p \in \R^2 \setminus 0 \},}\\
f_{A(x)}^{-1}(L') &=
\textstyle{\{ (\rho(|p|)p,-\frac{\pi-\psi(|p|)}{|p|}p) \suchthat p \in B_{2\epsilon}\}.}
\end{align*}
Here $\rho \in \smooth(\Rgeq,\R)$ is a cutoff function with $\rho(t) = 0$ for $t \leq \epsilon/4$
and $\rho(t) = 1$ for $t \geq \epsilon/2$. Note that $L'$ agrees with $\cotangent_xS^2 \cup S^2$ outside
$f_x(W)$, where $W = B_{2\epsilon} \times B_{\pi-\epsilon/4}$. The remaining portion of $L'$ can be written as
\begin{multline*}
f_x^{-1}(L') \cap W = \{ (y_1 \cos t, y_1 \sin t, y_2 \cos t, y_2 \sin t \suchthat \\
(y_1,y_2) \in C,\; t \in S^1 \} \cap W,
\end{multline*}
where $C \subset \R^2$ is the image of the embedding $c: \R^+ \longrightarrow \R^2$, $c(t) =
(\rho(\pi-\psi(t))t,-\psi(t))$. This is just the essential part of a Lagrangian handle
in $\R^4$.

Given two compact Lagrangian surfaces $L_1,L_2 \subset M$ which intersect transversely in a
single point $x$ and such that $L_2$ is a Lagrangian two-sphere, one can always find a
symplectic embedding $f: \disc \longrightarrow M$, for some $\epsilon>0$, such that
$f(S^2) = L_2$ and $f^{-1}(L_1) = \cotangent_x S^2 \cap \disc$ for some $x \in S^2$. The
argument above then proves that $\tau_{L_2}^{-1}(L_1)$ is Lagrangian isotopic to
$L_1 \surgery L_2$. \qed

The same argument shows that

\begin{proposition} \label{th:surgery-two}
Assume that $L_1$ is a Lagrangian sphere. Then $L_1 \surgery L_2$
is Lagrangian isotopic to $\tau_{L_1}(L_2)$. \qed \end{proposition}

\proof[Proof of Proposition \ref{th:braiding}]
Set $L_2' = \tau_{L_1}(L_2)$ and $L_1' = \tau_{L_2}^{-1}(L_1)$.
As a special case of Remark \ref{rem}\ref{item:naturality},
$\tau_{L_1} \tau_{L_2} \tau_{L_1}^{-1}$ is symplectically isotopic to $\tau_{L_2'}$.
Similarly $\tau_{L_2}^{-1} \tau_{L_1} \tau_{L_2}$ is symplectically isotopic to
$\tau_{L_1'}$. Propositions \ref{th:surgery-one} and \ref{th:surgery-two}
show that $L_1'$ and $L_2'$ are both Lagrangian isotopic to $L_1 \surgery L_2$
and hence Lagrangian isotopic to each other. It follows that $\tau_{L_1'}$
is symplectically isotopic to $\tau_{L_2'}$. \qed

\providecommand{\bysame}{\leavevmode\hbox to3em{\hrulefill}\thinspace}

\end{document}